\documentstyle[12pt]{article}
\newtheorem{theorem}{Theorem}[section]

\newtheorem{definition}[theorem]{Definition}
\newtheorem{remark}[theorem]{Remark}

\newenvironment{proof}{\noindent{\bf Proof.\/}}{$~\Box$}

\newcommand{\beq}{\begin{equation}}
\newcommand{\eeq}{\end{equation}}
\newcommand{\beqa}{\begin{eqnarray}}
\newcommand{\eeqa}{\end{eqnarray}}
\newcommand{\ba}{\begin{array}}
\newcommand{\ea}{\end{array}}

\def\dep#1#2{{\partial #1\over \partial #2}}

\def\depq#1#2{{\partial^2 #1 \over {\partial #2^2}}}
\def\depo#1{{\partial \over \partial #1}}
\def\into#1{{\int_{\Omega} #1}}         
         
\def\spazi{~~~~~~~~~~}

\begin{document}

\begin{center}
{\bf CONVERGENCE OF THE SPLITTING METHOD \\ 
FOR SHALLOW WATER EQUATIONS}
\footnote{Invited paper to the VII International Conference 
on {\it Artificial Intelligence and Information--Control Systems of Robots}, 
Institute of Computer Systems, Slovak Academy of Sciences, 
10--14 September, Bratislava (1997).} 
\end{center}

\vskip 0.5 truecm

\begin{center}
MARIA MORANDI CECCHI and LUCA SALASNICH
\vskip 0.5 truecm
Dipartimento di Matematica Pura ed Applicata \\
Universit\`a di Padova, Via Belzoni 7, I--35131 Padova, Italy\\
E-mail: mcecchi@math.unipd.it, salasnich@math.unipd.it
\end{center}

\vskip 1. truecm

\begin{center}
{\bf Abstract}
\end{center}

\vskip 0.5 truecm
\par
In this paper we analyze the 
convergence of the splitting method for shallow water equations. 
In particular we give an analytical estimation of the time step 
which is necessary for 
the convergence and then we study the behaviour 
of the motion of the shallow water in the Venice lagoon by using 
the splitting method with a finite element space discretization. 
The numerical calculations show that the splitting method 
is convergent if the time step of the first part is sufficiently small 
and that it gives a good agreement with the experimental data.

\section{Introduction}

The equations of a Newtonian (viscous) fluid are called 
{\it Navier--Stokes equations} 
and in the case of a inviscid fluid they give the Euler equations. 
The Navier--Stokes equations are given by
$$
{\partial {\bf u} \over \partial t}+ ({\bf u} \cdot {\bf \nabla}) {\bf u} + 
{\bf \nabla} p - \nu \nabla^2 {\bf u} = {\bf f} \; ,
$$
\beq
{\bf \nabla} \cdot {\bf u} = 0 \; ,
\eeq
where $p=P/\rho$ is the reduced pressure and 
$\nu =\mu/\rho$ is the reduced or kinematic viscosity$^{1}$. 
\par
The Navier--Stokes equations are to be integrated in the space--time domain 
$\Omega \times ]0,T[ \subset {\bf R}^3\times {\bf R}^+$ 
once an appropriate set of initial and boundary 
conditions has been defined. While the former conditions, in general, only 
need to reproduce a feasible shape of the current solution, the latter 
are to be set and treated with much care due to the effect they can have on 
the evolution process$^{2,3}$. 
\par
In the shallow water hypothesis one 
assumes that the characteristic horizontal 
scale $L$ for the motion (the wave length) 
is longer than the average height $\bar{h}$ of the fluid, 
i.e. $L>>\bar{h}$. In such hypothesis vertical 
acceleration and velocity 
are negligible and the flux becomes almost horizontal$^{4}$. 

\begin{definition}
The total height of the water is given by 
$$
h(x_1,x_2,t)= H(x_1,x_2)+\eta (x_1,x_2,t) \; ,
$$ 
where $H$ is the depth of the water in the stationary condition above 
the reference level $x_3=0$, and $\eta$ is the time--dependent difference, 
i.e. the height of the free surface. 
\end{definition}

\par
The external forces acting on the water are of extreme importance 
to determine the equations of motion of the system. 

\begin{definition}
The force of gravity is defined as 
$$
{\bf f}_g= \int_{\Omega} \; \rho \; {\bf g} \; d{\bf x} \; ,
$$ 
where ${\bf g}=[0,0,-g]^T$ is the vector acceleration of gravity. 
\end{definition}

If we consider large regions of water it is necessary 
to include other forces, like the Coriolis force and 
the Chezy force, which models the friction of the water 
at the bottom. 

\begin{definition}
The Coriolis force is defined as
$$
{\bf f}_{cor}=\int_{\Omega} \; {\bf \omega} \wedge {\bf u} \; d{\bf x} \; , 
$$
where $\omega =[0,0,w_3]^T$ is the rotation vector of the earth 
and $\omega_3=k_0$ is called Coriolis coefficient.  
\end{definition}

\begin{definition}
The Chezy force is defined as 
$$
{\bf f}_{ch}=\int_{\Omega} \; {g {\bf u} |u|\over k_1^2 h} d{\bf x} \; ,
$$
where $|u|=\sqrt{u_1^2+u_2^2}$, g is the scalar acceleration 
of gravity, $h$ is the total height of the water 
and $k_1$ is the Chezy coefficient. 
\end{definition}

\par
The first step to obtain the shallow water equations is to put 
$u_3$ and ${du_3/dt}$ equal to zero in the Navier Stokes equations. 
Then the third equation can be integrated between $-H$ and $\eta$. 
The system is yet 3--dimensional because 
the velocities $u_1$ and $u_2$ are functions also of the $x_3$ variable. 
It is possible to obtain a 2--dimensional system by performing 
the substitution
$$
u_1(x_1,x_2,x_3,t) \to a(\psi ) u_1(x_1,x_2,t) \; ,
$$
$$
u_2(x_1,x_2,x_3,t) \to a(\psi ) u_2(x_1,x_2,t) \; ,
$$
where 
$$
\psi = {x_3+H\over \eta + H} \spazi \hbox{ with } \spazi 
\int_0^1 a(\psi ) d\psi = 1 \; , 
$$
and then by integrating the new equations over the $x_3$ variable. 
In this way one can prove the following theorem (see Ref. 5 and 6 
for details). 

\begin{theorem}
Let us consider the inviscid shallow water 
with external forces given by the gravity force and the Coriolis force. 
The 2--dimensional viscid shallow water equations are 
\beq
{\partial {\bf U}\over \partial t} + {\partial {\bf F}_1 \over \partial x_1}  
+{\partial {\bf F}_2 \over \partial x_2} = {\bf R}_s
\label{eq:si1}
\eeq
where ${\bf U}=[\eta , u_1, u_2]^T$ is the vector of the conservative 
variables, the flow vectors ${\bf F}_1$ and ${\bf F}_2$ are 
\beq
{\bf F}_1 =[H u_1, g \eta , 0]^T \; , \;\;\;\;
{\bf F}_2 =[H u_2 , 0 , g \eta ]^T \; ,
\eeq
and the source vector ${\bf R}_s$ is given by 
\beq
{\bf R}_s = [ 0 , k_0 u_2 - {g u_1 |u|\over k_1^2 H} , 
- k_0 u_1 - {g u_2 |u|\over k_1^2 H} ]^T \; , 
\eeq
where $g$ is the acceleration of gravity, $k_0$ and $k_1$ are 
the coefficients of Coriolis and Chezy and $H(x_1,x_2)$ is the  
depth of the water in the stationary condition. 
\end{theorem}
\par
The shallow water equations are to be integrated in the space--time domain 
$\Omega \times ]0,T[ \subset {\bf R}^2\times {\bf R}^+$ 
once an appropriate set of initial and boundary 
conditions has been defined
\beq
{\bf U}(x_1,x_2,0)={\bf U}_0(x_1,x_2) \;\;\;\;\;\; 
\forall (x_1,x_2) \in \Omega \; ,
\eeq
\beq
{\bf U}(x_1,x_2,t)={\bf U}_{{\partial\Omega}}(x_1,x_2,t) \;\;\;\;\;\; 
\forall (x_1,x_2) \in {\partial\Omega} \; , \;
\forall t\in ]0,T[ \; .
\eeq 
While the initial conditions, in general, only 
need to reproduce a feasible shape of the current solution, the boundary 
conditions are to be set and treated with 
much care due to the effect they can have on 
the evolution process. 

\section{The splitting method}

A semi--implicit method is adopted for the solution of 
the system (\ref{eq:si1}), based on the following splitting 
(see Ref. 7 and 8)
\beq
{\bf F}_i = {\bf F}^*_i + {\bf F}_i^{**} \; , \;\;\;\;\; i=1,2
\eeq
and
\beq
{\bf U}^{(n+1)} = {\bf U}^{(n)} + \Delta {\bf U}^* 
+ \Delta {\bf U}^{**} \; ,
\eeq
where $\Delta {\bf U}^*$ and $\Delta {\bf U}^{**}$ are the  
increments of the solution vector. 
\par
In the iterative scheme we put ${\bf F}^*_i=0$ so that 
$F_i=F_i^{**}$, $i=1,2$, and the system (\ref{eq:si1}) can be divided in:
\beq
{\partial \Delta {\bf U}^* \over \partial t} = {\bf R}_s \; .
\label{eq:si2}
\eeq
and
\beq
{\partial \Delta {\bf U}^{**} \over \partial t}+ {\partial
{\bf F}^{**}_1 \over \partial x_1} + {\partial {\bf F}^{**}_2 \over
\partial x_2} ={\bf 0} \; .
\label{eq:si3}
\eeq
These equations are integrated in time, in turn, by using an explicit 
Taylor--Galerkin method$^{8}$ for (\ref{eq:si2}) and an implicit 
$\theta$--method$^{9}$ for (\ref{eq:si3}). 
\par
The equation (\ref{eq:si2}) is discretized in time by using 
a Taylor expansion to the second order
\beq
(\Delta {\bf U}^{*})^{(n+1)} =\tau \left(\dep {\Delta 
{\bf U}^*} t \right)^{(n)}+ {\tau^2 \over 2} \left(
\depq {\Delta {\bf U}^*} t \right) ^{(n)} \; ,
\label{eq:si4}
\eeq
where $\tau$ is the time step. From 
\beq
\dep {\Delta {\bf U}^*} t = {\bf R}_s \spazi \hbox{and} \spazi
\depq {\Delta {\bf U}^*} t =\dep {{\bf R}_s} t =
\dep {{\bf R}_s} {\Delta {\bf U}^*} \dep {\Delta {\bf U}^*} t \; ,
\label{eq:si5}
\eeq
the equation (\ref{eq:si4}) can be written
\beq
(\Delta {\bf U}^{*})^{(n+1)} =\tau ({\bf R}_s)^{(n)} +
{\tau^2 \over 2}({\bf G}{\bf R}_s)^{(n)}  \; ,
\label{eq:si6}
\eeq
where 
\beq
{\bf G}=\dep {{\bf R}_s} {\Delta {\bf U}^*} \; .
\eeq
Because of the computational complexity in the evaluation 
of the right term of equation (\ref{eq:si6}), we use a two--step version 
of the Taylor--Galerkin algorithm. This is given by
an approximation of $({\bf U})^{(n+1/2)}$ and ${{\bf R}_s}^{(n+1/2)}$
the Taylor expansion at step $(n+1/2)$
\beq
{\bf U}^{(n+1/2)}= {\bf U}^{(n)}+{\tau \over 2} 
({\bf R}_s)^{(n)} \; ,
\label{eq:si7}
\eeq
\beq
({\bf R}_s)^{(n+1/2)}=({\bf R}_s)^{(n)}+{\tau \over 2}
\left( \dep {{\bf R}_s} t\right)^{(n)}=
({\bf R}_s)^{(n)}+{\tau \over 2}({\bf G}{\bf R}_s)^{(n)} \; ,
\label{eq:si8}
\eeq
from which we obtain $({\bf G}{\bf R}_s)^{(n)}$ as
\beq
({\bf G}{\bf R}_s)^{(n)}={2\over \tau} \left[ ({\bf R}_s)^{(n+1/2)}-
({\bf R}_s)^{(n)} \right] \; .
\label{eq:si9}
\eeq
It follows that $(\Delta {\bf U}^{*})^{(n+1)}$ can be written as 
\beq
(\Delta {\bf U}^{*})^{(n+1)} =\tau ({\bf R}_s)^{(n+1/2)} \; ,
\eeq
and the explicit scheme of the splitting method results
\beq
{\bf U}^{(n+1)}={\bf U}^{(n)}+{\bf R}_s^{(n+1/2)} \; .
\eeq
\par
The equation (\ref{eq:si3}) is given, in explicit form, by
\begin{equation}
\begin{array}{l}
  \depo t \left( \Delta \eta^{**} \right) +
    \depo {x_1} \left( Hu_1 \right) +
    \depo {x_2} \left( Hu_2 \right) =0 \\
  \depo t \left( \Delta u_1^{**} \right) +
    \depo {x_1} \left( g\eta \right) =0 \\
  \depo t \left( \Delta u_2^{**} \right) +
    \depo {x_2} \left( g\eta \right) =0 \; .
\end{array}
\label{eq:si12}
\end{equation}
The time discretization is obtained by using the $\theta$-method 
\begin{equation}
\begin{array}{l}
  \left( \Delta \eta^{**} \right)^{(n+1)} +
  {\tilde{\tau}}\left[
     \depo {x_1} \left( Hu_1^{(n+\theta_1)} \right) +
     \depo {x_2} \left( Hu_2^{(n+\theta_1)} \right) \right] =0 \\
  \left( \Delta u_1^{**} \right)^{(n+1)} +
     {\tilde{\tau}}g \depo {x_1} \eta^{(n+\theta_2)} =0 \\
  \left( \Delta u_2^{**} \right)^{(n+1)} +
     {\tilde{\tau}}g \depo {x_2} \eta^{(n+\theta_2)} =0 \; 
\end{array}
\label{eq:si13}
\end{equation}
where ${\tilde \tau}$ is the time step, 
$\theta_1$ and $\theta_2$ are real parameters 
between $0$ and $1$ and 
\begin{equation}
\begin{array}{l}
  u_i^{(n+\theta_1)}=u_i^{(n)}+ \theta_1\left[
    \left( \Delta u_i^{*} \right)^{(n+1)} +
    \left( \Delta u_i^{**} \right)^{(n+1)} \right] \spazi i=1,2 \\
  \eta^{(n+\theta_2)}=\eta^{(n)}+\theta_2 (\Delta \eta^{**})^{(n+1)} \; .
\end{array}
\label{eq:si14}
\end{equation}
Note that the term $(\Delta \eta^{*})^{(n+1)}$ does not appear because the 
splitting method does not include variations in the water elevation 
due to (\ref{eq:si6}). Substituting (\ref{eq:si14}) into (\ref{eq:si13}), 
we can write 
\begin{equation}
\begin{array}{l}
  \left( \Delta \eta^{**} \right)^{(n+1)} +
     {\tilde{\tau}}\theta_1\sum_{i=1}^2 \depo {x_i}\left[
     H\left( \Delta u_1^{**} \right)^{(n+1)}  \right] = \\ \spazi \spazi
     -{\tilde{\tau}}\sum_{i=1}^2 \depo {x_i}\left( Hu_1^{(n)} \right) +
     {\tilde{\tau}}\theta_1\sum_{i=1}^2 \depo {x_i}\left[ 
     H\left( \Delta u_1^{*} \right)^{(n+1)}  \right]\\
  \left( \Delta u_1^{**} \right)^{(n+1)} +
     {\tilde{\tau}}\theta_2\depo {x_1} g \left[
     \left( \Delta \eta^{**} \right)^{(n+1)}\right] =
     -{\tilde{\tau}} g \depo {x_1} \eta^{(n)} \\
  \left( \Delta u_2^{**} \right)^{(n+1)} +
     {\tilde{\tau}}\theta_2 g \depo {x_2}\left[
     \left( \Delta \eta^{**} \right)^{(n+1)}\right] = 
     -{\tilde{\tau}}g \depo {x_2} \eta^{(n)} \; .
\end{array}
\label{eq:si15}
\end{equation}
Here $\left( \Delta u_1^{**} \right)^{(n+1)}$ and 
$\left( \Delta u_2^{**} \right)^{(n+1)}$ can be obtained from 
the second and third equation as a function of 
$\left( \Delta \eta^{**} \right)^{(n+1)}$ and then substituted 
in the first equation which becomes
\begin{equation}
\begin{array}{l}
  \left( \Delta \eta^{**} \right)^{(n+1)} -
     {\tilde{\tau}}^2\theta_1\theta_2g \sum_{i=1}^2 \depo {x_i}\left[
     H\depo {x_i}\left( \Delta \eta^{**} \right)^{(n+1)}  \right] = \\ \spazi
  -{\tilde{\tau}}\left\{\sum_{i=1}^2\depo {x_i}\left[ H\left( u_1^{(n)} +
     \theta_1\left( \Delta u_i^{*} \right)^{(n+1)}\right) \right] 
  -{\tilde{\tau}}\theta_1 g \sum_{i=1}^2 \depo {x_i}\left(
     H\depo {x_i} \eta^{(n)}\right) \right\} \; .
\end{array}
\label{eq:si16}
\end{equation}

\section{Time Convergence of the splitting method}

In this section we study analytically the time convergence of the 
the splitting method. Our study is motivated by the stability analysis 
of multilevel methods for the numerical simulation of turbulent flows  
performed by Roger Temam. He analyzed a simple model, namely a pair of 
ordinary coupled differential equations, by using numerical schemes 
based on different treatments and different time steps for the 
two variables of the system$^{10}$. 
\par
In our splitting method, starting from the time discretization 
of the previous section, we have 
$$
\eta^{(n+1)}= \eta^{(n)} + (\Delta \eta^{*})^{(n+1)} + 
(\Delta \eta^{**})^{(n+1)}\;
$$
\beq
u_1^{(n+1)}=u_1^{(n)} + (\Delta u_1^*)^{(n+1)} + (\Delta u_1^{**})^{(n+1)} \; 
\eeq
$$
u_2^{(n+1)}=u_2^{(n)} + (\Delta u_2^*)^{(n+1)} + (\Delta u_2^{**})^{(n+1)} \; .
$$
We observe that the first component of 
${\bf R}_s^{(n+1/2)}$ is zero so that $(\Delta \eta^*)^{(n+1)} = 0$. 
\par
It follows that the first part of the splitting scheme reads 
$$
\eta^{(n+1)}= \eta^{(n)} \;
$$
\beq
u_1^{(n+1)}=u_1^{(n)} + (\Delta u_1^*)^{(n+1)} \; 
\eeq
$$
u_2^{(n+1)}=u_2^{(n)} + (\Delta u_2^*)^{(n+1)} \; 
$$
where 
\beq
(\Delta u_1^*)^{(n+1)} = \alpha u_1^{(n)} + \beta  u_2^{(n)} \; 
\eeq
\beq
(\Delta u_2^*)^{(n+1)} = - \beta  u_1^{(n)} + \alpha u_2^{(n)} \; 
\eeq
with $\alpha = (-{\tau^2 k_0^2\over 2} - \tau  D + {\tau^2 D^2\over 2})$, 
$\beta = (\tau k_0 - \tau^2 D)$ and $D=g|u|/(k_1^2H)$. 
Here we suppose that $D$ is constant during the iterative process. 
The iteration scheme can be written as
\beq
\left( \begin{array}{c}  
      \eta^{(n+1)} \\ u_1^{(n+1)} \\ u_2^{(n+1)}
              \end{array}    \right) = J^* 
\left( \begin{array}{c}  
      \eta^{(n)} \\ u_1^{(n)} \\ u_2^{(n)}                
            \end{array}    \right) \; 
\eeq
where 
\beq
J^*= \left( \begin{array}{ccc} 
1 & 0          & 0           \\
0 & 1 + \alpha & \beta       \\
0 & - \beta    & 1 + \alpha  \\
\end{array}  \right) \; . 
\eeq

\begin{remark}
We observe that the first part of the splitting scheme and its 
matrix $J^*$ do not depend on $\theta_1$ and $\theta_2$. 
\end{remark}

Let us study the properties of the first part of the 
splitting scheme generated by the matrix $J^*$. 

\begin{definition}
The spectral radius of the matrix $J^*$ is given by 
$$
\rho (J^*) = max\{|\lambda_i|, \; i=1,2,3\} \; 
$$
where $\lambda_i$, $i=1,2,3$, are the eigenvalues of the matrix $J^*$. 
\end{definition}

\par 
By solving the eigenvalue equation we find that 
the eigenvalues of $J^*$ are $\lambda_1=0$ and 
$\lambda_{2,3}=(1+ \alpha ) \pm i \beta$.  

\begin{definition}
The iteration scheme induced by the matrix $J^*$ is called convergent 
if the spectral radius is such that $\rho (J^*)<1$. 
\end{definition}

Now we can prove the following theorem. 

\begin{theorem}
The first part of the splitting method, induced by the matrix $J^*$, is 
not convergent for $D=0$. For $D\neq 0$ it is convergent 
if and only if the following inequality is satisfied 
\beq
(k_0^4 + D^2 (4-k_0^2)+4D^2) \tau^3 
- 4D(D^2 +2k_0 -k_0^2) \tau^2 + 8D^2 \tau - 8 D < 0 \; , 
\eeq
where $D=g|u|/(k_1^2H)$. 
\end{theorem}

\begin{proof}
The convergence condition for the iterative process 
is such that the spectral radius of $J^*$ is less than one. 
We have seen that the eigenvalues of $J^*$ are $\lambda_1=0$ 
and $\lambda_{2,3}=(1+\alpha)\pm i \beta$. It follows that the convergence 
condition is $(1+\alpha )^2 + \beta^2 < 1$, where $\alpha=
(\tau^2 {k_0^2\over 2} - \tau  D + \tau^2 {D^2\over 2})$ 
and $\beta=(\tau k_0 - \tau^2 D)$, and 
this condition gives the inequality of the theorem. 
In particular, if $D=0$ we have the inequality $k_0^4 < 0$ 
which is never satisfied. 
\end{proof}

\begin{remark}
The convergence of the iteration scheme generated by the matrix $J^*$ 
does not depend on ${\tilde \tau}$, $\theta_1$ and $\theta_2$ 
because these parameters do not appear in the eigenvalues 
of $J^*$ and in the inequality of the theorem 3.3. 
\end{remark}

\par
The equation associated to the inequality of the theorem 3.3 is given by 
\beq 
a \tau^3 - b \tau^2 + c \tau - d = 0 \; , 
\eeq
where $a=(k_0^4 + D^2 (4-k_0^2)+4D^2)$, 
$b=4D(D^2 +2k_0 -k_0^2)$, $c=8D^2$ and $d=8D$. 
This equation has two complex conjugate solutions and 
a real solution given by 
\beq 
\tau_c = {b\over 3 a} - {2^{1/3} (-b^2 + 3 ac)\over 3 a q^{1/3}} 
+ {1\over 2^{1/3} 3 a} q^{1/3} \; , 
\eeq
where $q=2b^3-9abc+27a^2d+\sqrt{4(-b^2+3ac)^3+(2b^3-9abc+27a^2d)^2}$. 
It follows that the iteration scheme induced by the matrix $J^*$ is 
convergent if and only if the following inequality is satisfied 
\beq
\tau < \tau_c \; .
\eeq
In particular, with the positions 
$g = 9.81$ m/sec$^2$, $k_0 = 10^{-4}$ sec$^{-1}$, 
$k_1 = 40$ m$^{1/2}$sec$^{-1}$, $|u| = 0.1$ m/sec and $H=0.1$ m, 
we find $\tau < \tau_c = 5.41$ sec. \\
We shall use this estimation 
for the numerical implementation of the splitting algorithm. 
\par
In conclusion, we see that the convergence of the first part 
of the splitting method depends strongly 
on the time step $\tau$. This time step must be very small 
to ensure the convergence of the iteration scheme. Moreover 
the inclusion of the Chezy force is essential for the convergence. 
In fact, if the Chezy force is not included ($D=0)$ then the iterative 
process diverges. 
\par
We observe that the problem of spatial discretization 
can be studied with the Finite Element Method (see next section). 
As shown by Ciarlet and Raviart$^{11}$, 
the interpolating function of Finite Element Method 
controls completely the spatial convergence. It follows that, 
to analyze the time convergence of the second part of the splitting 
method, the spatial dependence can be neglected. 
In this way, we get $(\Delta u_i^{**})^{(n+1)}=0$, $i=1,2$, and 
the splitting scheme reads 
$$
\eta^{(n+1)}= \eta^{(n)} + (\Delta \eta^{**})^{(n+1)} \;
$$
\beq
u_1^{(n+1)}=u_1^{(n)} + (\Delta u_1^*)^{(n+1)} \; 
\eeq
$$
u_2^{(n+1)}=u_2^{(n)} + (\Delta u_2^*)^{(n+1)} \; 
$$
where 
$$
(\Delta \eta^{**})^{(n+1)} = - {\tilde{\tau}} 
[{\partial H\over \partial x_1} 
+{\partial H\over \partial x_1} 
\theta_1 \alpha 
 - {\partial H\over \partial x_2} 
\theta_1 \beta ] u_1^{(n)} - {\tilde{\tau}} [{\partial H\over \partial x_1} 
\theta_1 \beta + {\partial H\over \partial x_2} 
+ {\partial H\over \partial x_2} \theta_1 \alpha ] u_2^{(n)} \; 
$$
\beq
(\Delta u_1^*)^{(n+1)} = \alpha u_1^{(n)} + \beta  u_2^{(n)} \; 
\eeq
$$
(\Delta u_2^*)^{(n+1)} = - \beta  u_1^{(n)} + \alpha u_2^{(n)} \; 
$$
with $\alpha = (-{\tau^2 k_0^2\over 2} - \tau  D + {\tau^2 D^2\over 2})$, 
$\beta = (\tau k_0 - \tau^2 D)$ and $D=g|u|/(k_1^2H)$. 
It follows that the iteration scheme can be written as
\beq
\left( \begin{array}{c}  
      \eta^{(n+1)} \\ u_1^{(n+1)} \\ u_2^{(n+1)}
              \end{array}    \right) = J^{**} 
\left( \begin{array}{c}  
      \eta^{(n)} \\ u_1^{(n)} \\ u_2^{(n)}                
            \end{array}    \right) \; 
\eeq
where 
\beq
J^{**}= \left( \begin{array}{ccc} 
 1 & - {\tilde{\tau}} [{\partial H\over \partial x_1} +
{\partial H\over \partial x_2}\theta_1 \alpha 
- {\partial H\over \partial x_2} \theta_1 \beta ] &  
- {\tilde{\tau}} [{\partial H\over \partial x_2} + 
{\partial H\over \partial x_1}\theta_1 \beta 
- {\partial H\over \partial x_2} 
\theta_1 \alpha ]     \\
0 & 1 + \alpha & \beta \\
0 & - \beta   & 1 + \alpha \\
\end{array}  \right) \; . 
\eeq

\begin{remark}
The iteration scheme induced by the matrix $J^{**}$ 
does not depend on $\theta_2$. 
\end{remark}

\par
It is easy to see that the eigenvalues of $J^{**}$ are the same of $J^*$. 
They are $\lambda_1=0$ and $\lambda_{2,3}=(1+ \alpha ) \pm i \beta$.  
It means that also the convergence of the iteration scheme 
generated by the matrix $J^{**}$ 
does not depend on ${\tilde \tau}$, $\theta_1$ and $\theta_2$ 
because these parameters do not not appear in the eigenvalues 
of $J^{**}$. In conclusion, we do not have limitations for ${\tilde \tau}$. 
In practice, because of the effect of the error propagation due to the 
numerical approximations, we use ${\tilde \tau}\simeq 5$ minutes.  

\section{Finite element space discretization}

The spatial approximation with the finite element method is 
obtained by using linear form functions for the integer step 
$(n,n+1,\dots)$ and constant functions for the half step 
$(n-1/2,n+1/2,\dots)$. In this way the equation  
(\ref{eq:si6}) reads (see also Ref. 12 and 13)
\begin{eqnarray}
({\bf M}\Delta {\bf U}^{*})_j^{(n+1)} = \tau\into {({\bf R}_s)^{(n)}
\phi_j d{\bf x} } +{\tau^2 \over 2} \into {{2\over \tau} \left[
({\bf R}_s)^{(n+1/2)}-({\overline {\bf R}}_s)^{(n)} \right] \phi_j 
d{\bf x} }
\label{eq:si10}\\
=\tau\into{
    \left\{ 
      ({\bf R}_s)^{(n+1/2)} +
      \left[ ({\bf R}_s)^{(n)} -({\overline {\bf R}}_s)^{(n)}
      \right] 
    \right\} 
    \phi_j d{\bf x} 
} \spazi
j=1,2,\cdots ,N\nonumber 
\end{eqnarray}
where $N$ is the total number of nodes, $\phi_j$ is the weight linear function 
of the node $j$, the bar denotes mean values calculated on the element and 
\begin{equation}
{\bf M}=[M]_{ij}=
\into{\phi_i \phi_j d{\bf x} }  \spazi i,j=1,2,\cdots ,N
\label{eq:si11}
\end{equation}
is the mass matrix.
\par
By using again linear triangular elements for the spatial discretization 
with the finite elements method, and by using the Green formulas 
to the terms which include the second derivatives, the variational 
formulation (\ref{eq:si16}) gives the following system of discrete 
equations 
\begin{equation}
\begin{array}{l}
  \left( {\bf M}+{\tilde{\tau}}^2 g\theta_1\theta_2{\bf S}\right)
  (\Delta p^{**})^{(n+1)} = \\ \spazi-{\tilde{\tau}}
\left\{\sum_{i=1}^2{\bf Q}_i
  \left[ H\left( u_i^{(n)}+\theta_1\left( \Delta u_i^{*} \right)^{(n+1)}
  \right)\right]+{\tilde{\tau}}\theta_1{\bf S}p^{(n)}\right\}
\end{array}
\label{eq:si17}
\end{equation}
where
\begin{equation}
  {\bf S}=\sum_{i=1}^2\left(\into{ 
{\partial [\phi] \over \partial x_i} H {\partial [\phi]^T \over \partial x_i} 
 d{\bf x} }\right)\spazi \hbox{and} \spazi
  {\bf Q}_i = \into [\phi] {\partial [\phi]^T \over \partial x_i} 
d{\bf x} \; .
\label{eq:si18}
\end{equation}
When $\left( \Delta \eta^{**} \right)^{(n+1)}$ is evaluated by 
(\ref{eq:si17}), the second and third equations of (\ref{eq:si15})
can be solved in 
$\left( \Delta u_1^{**} \right)^{(n+1)}$ and 
$\left( \Delta u_2^{**} \right)^{(n+1)}$; 
their discretization to finite elements is given by 
\begin{equation}
  {\bf M}(\Delta u_i^{**})^{(n+1)}=-{\tilde{\tau}} g {\bf Q}_i
\left[ \eta^{(n)}+\theta_2( \Delta \eta^{**})^{(n+1)}\right]\spazi i=1,2 \; .
\label{eq:si19}
\end{equation}
It is easy to see that this set of equations is similar to 
(\ref{eq:si10}); as a consequence it can be solved in the same way.

In summary, the solution at each time step implies the following 
operations:
\begin{description}
\item [a)]~~solve (\ref{eq:si10}) for $( \Delta u^{*})^{(n+1)}$ and  
     $( \Delta u^{*})^{(n+1)}$(here $( \Delta \eta^{*})^{(n+1)}=0$ 
     because the first component of ${\bf R}_s$ is zero);
\item [b)]~~calculate $( \Delta \eta^{**})^{(n+1)}$ from (\ref{eq:si17});
\item [c)]~~solve (\ref{eq:si19}) for $( \Delta u_1^{**})^{(n+1)}$ and 
     $( \Delta u_2^{**})^{(n+1)}$;
\item [d)]~~calculate ${\bf U}^{(n+1)} = {\bf U}^{(n)} +
     (\Delta {\bf U}^*)^{(n+1)}+(\Delta {\bf U}^{**})^{(n+1)}$.
\end{description}

\section{Numerical calculations}

In this section we describe our calculations of the motion 
of the water into the Venice lagoon$^{12,13}$. 
The equation considered are the shallow water equations.  
All simplification are obtained considering different depth 
of the lagoon and in particular a finer discretization is 
considered where the water runs with greater speed 
namely along the more important canals. 
The numerical method is based on the splitting scheme and 
the finite element discretization presented in the previous sections. 
\par
To analyze the Venice Lagoon we must add the wind effect. Let 
$v$ the wind velocity respect the water and $\xi$ an adimensional constant, 
the wind term is given by $\xi |v|v/H$. 
In the approximation that the water velocity is negligible respect 
the wind velocity, the wind term is an external parameter to the system, i.e. 
the wind velocity respect the soil. We take a time--dependent 
wind velocity but space--independent$^{14}$. 
\par 
Following our analytical estimation of the time step which is 
necessary for the convergence of the splitting scheme, 
we choose a time step $\tau = 3$ seconds for the explicit part and 
${\tilde{\tau}} =5$ minutes for the implicit part. 
In this way the model is very stable and capable to preserve 
the water quantity of the system. 
We decompose the Venice lagoon into 10 
subdomains and the nodes of the interfaces collected into a unique set. 
The number of nodes of the entire discretization is 1967 and the elements are 
3423. The results of the calculation obtained using the domain decomposition 
method without parallelization have approximately a speed--up of 25 per cent 
over a Conjugate Gradient solver used 
on a not--decomposed domain$^{12,13,17}$. 
\par
The model has been verified by using the experimental data 
of the Ufficio Idrografico and Mareografico di Venezia 
(month of Fabruary 1994)$^{15}$.
The data have been introduced at the mouth of the lagoon 
as boundary conditions on the open boundary. Therefore a comparison 
has been made between the simulation and the real data on the 
internal measurement stations of the lagoon$^{16}$. 
The experimental and simulated data present a good agreement 
in all the internal measurement stations$^{12,13,17}$. 
 
\section*{References}
\parindent=0. pt

1. A. J. Chorin and J. E. Marsden, 
{\it A Mathematical Introduction to Fluid Mechanics} 
(Springer, Berlin, 1980).

2. L. Quartapelle, 
{\it Numerical Solution of the Incompressible 
Navier--Stokes Equations} (Birkhouse, Bessel, 1993).

3. R. Temam, {\it Navier-Stokes Equations} 
(North Holland, Amsterdam, 1977).     

4. G.F. Carey and J.T. Oden, {\it Finite Elements: Fluid Mechanics} 
(Prentice--Hall, New Jersey, 1986).        

5. J. Peraire, O.C. Zienkiewicz and K. Morgan, 
"Shallow water problems: a general explicit formulation", 
Int. J. Numer. Meth. Eng. {\bf 22}, 547--574 (1986).

6. M. Morandi Cecchi, "Study of the behaviour of oscillatory waves 
in a lagoon", Int. J. Numer. Meth. Eng. {\bf 27}, 103--112 (1989). 

7. O.C. Zienkiewicz and P. Ortiz, "A split--characteristic based 
finite element model for shallow water equations", 
Int. J. Numer. Methods Fluids {\bf 20}, 1061--1080 (1995).

8. O.C. Zienkiewicz and P. Ortiz, "An Improved Finite Element Model 
for Shallow Water Problems", 
in {\it Finite Element Modelling of Enviromental Problems}, 
ch. 4, pp. 61--84 (Wiley, New York, 1995).

9. R. Glowinski, "Finite element methods for the numerical 
simulation of incompressible viscous flow. Introduction to the 
control of Navier--Stokes equations", 
Lect. Appl. Math. {\bf 28}, 219--301 (1991).

10. R. Temam, "Multilevel Methods for the Simulation of Turbulence", 
J. Comput. Phys. {\bf 127}, 309--315 (1996).

11. P.G. Ciarlet and P.A. Raviart, "General Lagrange and Hermite 
interpolation in $R^n$ with applications to finite element methods", 
Arch. Rational Mech. Anal. {\bf 46}, 177--199 (1972).

12. M. Morandi Cecchi, A. Pica and E. Secco, "Tidal 
flow analysis with core minimization", 
in the 8th International Conference on Finite Elements in Fluids: 
Ed. by K. Morgan {\it et al.}, pp. 990--998 (Pineridge Press, U.K., 1993).

13. M. Morandi Cecchi, A. Pica and E. Secco, "Domain decomposition 
of a finite element model of a lagoon", 
in the 8th International Conference on Finite Elements in Fluids: 
Ed. by K. Morgan {\it et al.}, pp. 1016--1036 (Pineridge Press, U.K., 1993).

14. G. Di Silvio, "Modelli matematici per lo studio della 
propagazione di onde lunghe e del trasporto di materia nei corsi 
d'acqua e nelle zone costiere", 
in Equazioni Differenziali dell'Idrologia e dell'Idraulica, 
Ed. by G. Casadei {\it et al.}, pp. 11--41 (Patron, Bologna, 1979).

15. L. D'Alpaos and G. Di Silvio, 
{\it Le correnti di marea nella Laguna di Venezia}, 
(Ministero dei Lavori Pubblici, 1979).    

16. C.N.R. Grandi Masse, 
{\it Foglio accompagnatorio dei rilevamenti digitalizzati delle profondit\`a
nella laguna di Venezia} (C.N.R. Grandi Masse, Venezia, 1970).   

17. M. Morandi Cecchi and L. Salasnich, 
"Shallow water theory and its application to the Venice lagoon", 
Preprint 3/97, Dep. of Pure and Appl. Math., Univ. of Padova, 
to be published in the Proceedings of the Symposium on 
Advances in Computational Mechanics, 13--15 January 1997, 
Austin (Texas). 

\end{document}